# TEXTURE SYNTHESIS AND NONPARAMETRIC RESAMPLING OF RANDOM FIELDS

By Elizaveta Levina and Peter J. Bickel

*University of Michigan and University of California, Berkeley*

This paper introduces a nonparametric algorithm for bootstrapping a stationary random field and proves certain consistency properties of the algorithm for the case of mixing random fields. The motivation for this paper comes from relating a heuristic texture synthesis algorithm popular in computer vision to general nonparametric bootstrapping of stationary random fields. We give a formal resampling scheme for the heuristic texture algorithm and prove that it produces a consistent estimate of the joint distribution of pixels in a window of certain size under mixing and regularity conditions on the random field. The joint distribution of pixels is the quantity of interest here because theories of human perception of texture suggest that two textures with the same joint distribution of pixel values in a suitably chosen window will appear similar to a human. Thus we provide theoretical justification for an algorithm that has already been very successful in practice, and suggest an explanation for its perceptually good results.

**1. Introduction.** Texture is one of the central concepts in computer vision and image analysis. The term is generally used to refer to images of repeated patterns with local variations, such as waves, sand or human tissue. The stochastic nature of texture variations, not necessarily present in other real images, makes it a particularly natural area for applying statistical methods. While many texture algorithms are deterministic and based on heuristics rather than probability models, we found that a statistical framework can help to understand, justify and improve such algorithms; in turn, issues arising in texture algorithms can lead to questions of general statistical interest.

The problem of texture synthesis, which lies at the intersection of computer vision and computer graphics, is the problem of producing a new texture image which looks like a given texture, but is not exactly the same. It is frequently used in computer graphics to "paint" textures on surfaces, and can also be used for image compression, where the whole texture can be recreated from a small sample. The point of departure for our research is a simple and very popular heuristic resampling algorithm for texture synthesis [9] which produces excellent visual results but has no theoretical justification or statistical framework.

We formalize this algorithm in the framework of resampling from random fields and prove that it provides consistent estimates of the joint distribution of pixels in a window of specified size. The interest in the joint distribution of pixels in a window comes from theories of human perception of texture. The study of human pre-attentive texture discrimination was pioneered by Julesz in the 1960s and 1970s [12, 13, 14]. His original conjecture was that textures appear indistinguishable to humans if they have identical first- and second-order statistics, and was later extended to higher-order statistics (i.e., joint distributions of pairs, triples, etc.). When textures are viewed as random fields on a lattice, they are often assumed stationary and Markovian, in which case the distribution of $k$ pixels in the Markov neighborhood determines the $k$-order statistics, and the whole distribution.

A more modern view of texture perception is that the cells in the visual cortex respond to primitive stimuli like bars, edges, dots, and so on, at different frequencies and orientations. Psychophysical and neurophysiological experiments suggest that the brain performs multichannel spatial frequency and orientation initial analysis of any image formed on the retina and not just texture [6, 11]. These and other similar findings inspired the multichannel filtering approaches, which use distributions of filter responses for texture discrimination, with texture boundaries corresponding to sudden changes in the intensity of "firing" of some of the filters. A comprehensive texture perception model based on this idea was proposed by Malik and Perona [19], and many filter-based methods were developed subsequently. This view also supports the claim that the joint distribution of $k$ neighboring pixels determines texture appearance, since the joint distribution of pixel intensities in the filter support window determines the distribution of filter responses. In our view, these two interpretations of human perception complement each other, and both point to the joint distribution of pixels in a window as the key quantity.

In computer vision, texture synthesis algorithms are ultimately evaluated by human visual assessment of synthesized texture. Here we provide a proof that, according to theories of human perception, the algorithm of Efros and Leung can be expected to produce good visual results. To the best of our knowledge, the only other texture synthesis algorithm in the literature with a mathematical justification is FRAME [28], but, unfortunately, it does not produce very good visual results in practice, whereas the algorithm considered here does.

This paper converts the Efros and Leung algorithm into a formal bootstrap scheme for resampling stationary random fields. The bootstrap techniques for stationary random fields in the statistical literature are primarily used for estimating the mean and the variance of the random field, a goal very different from synthesis or estimating the joint distribution. The main tool used in this context is the moving block bootstrap (MBB) and its variants. MBB was first introduced for time series [15, 18] and extended to general random fields by Politis and Romano [21]. It is based on resampling blocks independently and concatenating them, rather than resampling by conditioning on the neighboring blocks, which is the main difference between our bootstrap algorithm and MBB. For time series, bootstrapping by conditioning on the past has been introduced by Rajarshi [24] and Paparoditis and Politis [20]; here we extend their methods to stationary random fields.

This paper is organized as follows. In Section 2 we give some background on texture synthesis and introduce the algorithm of Efros and Leung [9]. In Section 3 we formalize the algorithm in the framework of resampling from stationary random fields, and introduce a special case of Markov mesh models, which motivate a natural ordering on the plane. In Section 4 we show that both the Markov mesh version and the original algorithm produce consistent nonparametric estimates of the joint distribution of pixels in a patch, though the patch sizes differ for the two algorithms. This result is proved under the assumptions that the texture is a sample from a stationary mixing random field with a smooth density with compact support, and some minor regularity conditions. Section 5 concludes with discussion, and the Appendix contains all the proofs.

**2. The nonparametric sampling algorithm and previous work in texture synthesis.** There has been a surge of interest in texture synthesis in the past decade, when advances in computing allowed using many computationally intensive algorithms that could not have been implemented before. The many different methods of texture synthesis can be broadly divided into three categories. The first and oldest group of methods is model-based, with the main modeling tool being Markov random fields (MRF's) [2, 4]. In the early MRF work only a few parameters could be fitted because of computational difficulties, and those models usually did not capture the complexity of real textures. As the number of parameters increases, the synthesized textures begin to look more realistic, but it also becomes hard to estimate the parameters reliably.

The other broad category of texture synthesis methods is based on feature matching. Typically, these methods start from a white noise image and force it to match some set of statistics of the original texture image, such as distributions of filter responses [5, 10, 23, 25]. Feature matching methods tend to work well on stochastic textures but have difficulties with highly structured textures. Another difficulty is that they typically require some number of iterations to converge but iterating too many times leads to deterioration of the synthesized image.

There are some methods that use both MRF models and feature matching, such as the FRAME model by Zhu, Wu and Mumford [28]. It provides a solid theoretical base for combining MRF's with feature matching, and Wu, Zhu and Liu [27] showed that FRAME is the natural way to establish equivalence between these two approaches. However, its visual results on real textures are unfortunately far from perfect.

A new class of heuristic methods of texture synthesis has been developed over the past few years, started by the algorithm of Efros and Leung [9]. Many variations of their method have been published that speed up and optimize the original algorithm in different ways [8, 17, 26]. In all these works, however, the basic resampling principle of Efros and Leung [9] remains unchanged, and even the original version has been very successful on a wider range of textures than any of the previous methods.

The Efros and Leung algorithm is based on resampling from the random field directly, without constructing an explicit model for the distribution. It is motivated by an MRF model, that is, by the idea that the value of a given pixel only depends on the values of its neighbors, though it is not explicitly assumed that the underlying texture distribution is an MRF.

The algorithm starts with a random "seed" from the original image, typically a small square patch, and proceeds to grow the image from the seed outward, layer by layer, spiraling around and adding one pixel at a time. To synthesize pixel $X$, one conditions on $O(X)$, the part of the Markov neighborhood of $X$ (taken to be $w \times w$ square) that has been filled in before $X$ (see Figure 1). The conditional distribution of $X$ given $O(X)$ is never constructed explicitly. Instead, it is resampled directly in the following way: for all pixels $Y_i$ in the observed image compute the distance $d(O(X), O(Y_i))$, for all neighborhoods $O(Y_i)$ of the same size and shape as $O(X)$. The distance is measured by the sum of squared differences between pixel intensities, weighted by a Gaussian weight function to emphasize the importance of close neighbors. Let

$$d_0 = \min_i d(O(X), O(Y_i))$$

be the distance to the best match in the observed image. Define the set of "good" matches by

$$S = \{Y : d(O(X), O(Y)) \leq (1 + \varepsilon) d_0\}.$$

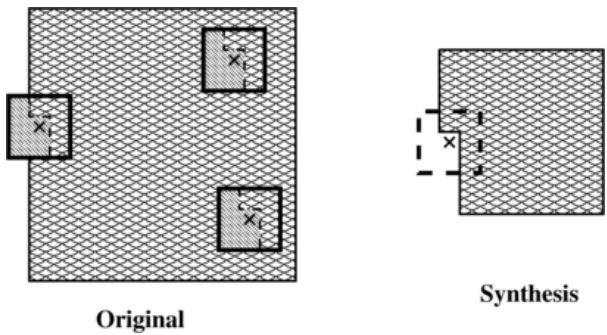

Fig. 1. *The nonparametric resampling algorithm.*

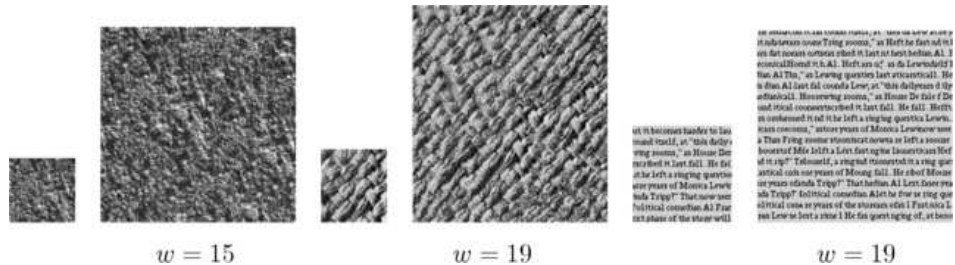

FIG. 2. *Some synthesis examples. The original (smaller) texture sizes are* $151 \times 148$, $54 \times 60$ *and* $113 \times 110$ *pixels, respectively.*

Finally, select the value for $X$ uniformly from pixel values in $S$.

Here $\varepsilon$ is a tuning parameter set by Efros and Leung to be $\varepsilon = 0.1$ (presumably by trial and error), and it is not meant to be changed by the user; this value was used in all the Efros and Leung results shown below.

This algorithm is very simple to implement, and can be used to synthesize any size or shape of the desired texture, or fill holes in an existing texture. It has worked well on both stochastic and structured textures (see Figure 2 for some examples). Note that highly structured textures require larger window sizes than more stochastic textures, and in general, the success of the algorithm depends on the neighborhood window being big enough to capture the local structure correctly, as shown in Figure 3 [9]. The smaller the window, the more "stochastic" the synthesized image will appear. This issue is discussed further in Section 5; choosing the window automatically is beyond the scope of this paper.

Although the algorithm produced impressive results on a large number of various textures, it also produced a few failures, discussed by Efros and Leung [9]. It appears to fail when it gets into a part of the search space with no good matches; in that case, it starts sampling randomly and produces texture that looks rather like white noise, but the chance of that happening is small. For most practical purposes, the algorithm works quite well, and, particularly with later computational speed-ups, is the current state of the art in texture synthesis.

**3. Formalizing the resampling scheme.** In this section we set up a formal bootstrap scheme along the lines of the synthesis algorithm heuristic. Our scheme is an extension to random fields of a $p$-order Markov bootstrap algorithm for time series by Paparoditis and Politis [20], and we use many of their techniques in the proofs. Their algorithm is, in turn, an extension of a first-order Markov bootstrap of Rajarshi [24]. As in the texture synthesis algorithm, the Markov assumption on the original time series is not needed in Paparoditis and Politis [20], but the bootstrapped time series reproduce the $p$-order dependence structure accurately.

For the case of random fields, we will first consider a Markov mesh model (MMM), a special case of MRF, which, unlike a general MRF on the plane, has a natural notion of the past.

3.1. *The resampling algorithm for Markov mesh models.* MMM's (also known as Picard random fields) were introduced by Abend, Harley and Kanal [1] and have been used for a variety of applications. In particular, Popat and Picard [22] used a parametric MMM model for texture synthesis, and so did Cressie and Davidson [3]. In both cases, however, results for natural textures were of low quality due to the small size of the conditioning neighborhood. Fitting all the parameters required for a larger neighborhood was computationally infeasible at the time, and the accuracy of estimating so many parameters would have been low in any case.

To define MMM's, let $\{X_t, t \in [1, \infty)^2\}$ be a real-valued random field. For a point $t = (t_1, t_2)$, define the index set

$$U_t = \{u : \max(1, t_1 - w + 1) \leq u_1 \leq t_1, \max(1, t_2 - w + 1) \leq u_2 \leq t_2, u \neq t\}$$

to be a square of size $w \times w$ with $t$ as the bottom right corner, $t$ itself excluded; notice that for the first $w - 1$ rows and columns $U_t$ has to be

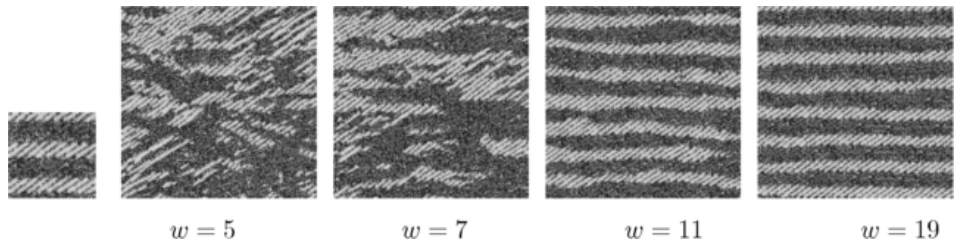

FIG. 3. *Synthesis results with different window sizes. The original image is $73 \times 71$ pixels.*

truncated. Let

$$W_t = \{u : 1 \le u_1 < t_1\} \cup \{u : 1 \le u_2 < t_2\}$$

be everything to the left or above $t$ (see Figure 4). Then a Markov mesh model assumes

$$P(X_t | X_{W_t}) = P(X_t | X_{U_t}).$$

MMM's are a special case of Markov random fields and here the corresponding Markov neighborhood $N_t$ is a $(2w - 1) \times (2w - 1)$ square centered at $t$, that is, $P(X_t | X_{-t}) = P(X_t | X_{N_t})$.

If the texture synthesis algorithm is to be motivated by a MMM, the natural way to fill in the pixels is to start in the upper left corner and proceed in raster order, filling in row by row. Suppose we observe the MMM field $X_t$ on the index set $[1, T_1] \times [1, T_2]$. Let $U_t(s)$ be the index set $U_t$ shifted so that its bottom right corner is $s$: $U_t(s) = (U_t - t + s)$. For convenience, define the $p$-dimensional vectors $\mathbf{Y}_t = X_{U_t}$ and $\mathbf{Y}_t(s) = X_{U_t(s)}$. Stationarity implies that $\mathbf{Y}_t(s)$ are identically distributed for all $t$ and $s$.

There are $w^2 - 1$ possible shapes of $U_t$ (various truncations of the $w \times w$ square are needed at the boundaries). For each shape consisting of $p$ components ($1 \le p \le w^2 - 1$), let $W^{(p)}$ be a kernel on $\mathbb{R}^p$. The kernel can be scaled by a resampling width $b$, $W_b^{(p)}(\mathbf{y}) = b^{-p} W^{(p)}(\mathbf{y}/b)$, and satisfies some general smoothness assumptions we state in Section 4.1. In the synthesis examples below, we use the Gaussian kernel $W^{(p)}(\mathbf{y}) = (2\pi)^{-p/2} \exp(-\|\mathbf{y}\|^2/2)$. Now we have all the components to proceed to

THE MMM RESAMPLING ALGORITHM.

1. Select a starting value for $\{X_t^* : t_1 \le w, t_2 \le w\}$, the top left $w \times w$ square. Typically the starting value will be a $w \times w$ square selected from the observed field $X_t$ at random.
2. Suppose $X_t^*$ has been generated for $\{t : t_1 < u\} \cup \{t : t_1 = u, t_2 < v\}$, that is, $u - 1$ rows are filled in completely, and row $u$ is filled up to column

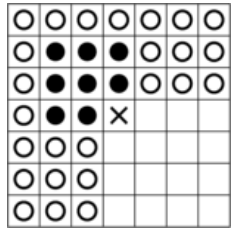

FIG. 4. *Conditional independence structure in the Markov mesh model.*

   $v$. To generate the next value $X_t^* = X_{(u,v)}^*$, let $N$ be a discrete random variable with probability mass function

$$P(N = s) = \frac{1}{Z} W_b^{(p)}(\mathbf{Y}_t^* - \mathbf{Y}_t(s)),$$

where $Z = \sum_s W_b^{(p)}(\mathbf{Y}_t^* - \mathbf{Y}_t(s))$ is a normalizing constant, $p = |\mathbf{Y}_t^*|$ is the size of the "past" of $X_t^*$, and $s$ ranges over all values $s$ such that $U_t(s) \subset [1, T_1] \times [1, T_2]$, that is, all locations where the conditioning neighborhood fits within the observed texture field.

3. Let $X_{(u,v)}^* = X_N$.

3.2. *Formalizing the general algorithm.* The MMM version of the algorithm contains two modifications of the original algorithm of Efros and Leung [9]: the order in which the pixels are filled in the synthesized texture (raster instead of spiral), and the weights with which the pixels are resampled (kernel weights instead of uniform sampling from all matches within $\varepsilon$). A number of comparisons we give in Section 3.3 show that both versions produce reasonable and fairly similar results; however, the spiral ordering of the original algorithm tends to have fewer problems with error propagation and produces somewhat more visually pleasing pictures. Therefore it is of interest to investigate the consistency properties of the spiral algorithm as well.

Here we will think of texture as a stationary random field on $\mathbb{Z}^2$ rather than $\mathbb{N}^2$. Let us order all locations $t \in \mathbb{Z}^2$ in the spiral order $t_0 \prec t_1 \prec \cdots$ starting at the origin and going around clockwise:

$$t_0 = (0,0), \quad t_1 = (1,0), \quad t_2 = (1,-1), \quad t_3 = (0,-1),$$
$$t_4 = (-1,-1), \quad t_5 = (-1,0), \quad t_6 = (-1,1), \quad t_7 = (0,1), \ldots.$$

To avoid centering problems, we will only look at conditioning on windows with an odd number of pixels along the side of the square, $(2w-1) \times (2w-1)$. The first $m^2$ pixels $(t_0, t_1, \ldots, t_{m^2-1})$ will be filled in by the seed, say also of size $m = 2w - 1$, and the subsequent pixels will be filled in one by one according to the spiral ordering. Apart from the ordering, the resampling scheme is exactly the same as for the MMM algorithm.

For all $t \in \mathbb{Z}^2$, let $X_t$ be the pixel intensity at location $t$. Let

$$\|t\|_\infty = \max(|t_1|, |t_2|)$$

be the $l_\infty$ norm on the plane. Let

$$U_t = \{s : \|t - s\|_\infty < w, s \prec t\}$$

be the part of $(2w-1) \times (2w-1)$ window with $t$ in the center that is filled in before $t$. Let $U_t(s)$ be the index set $U_t$ shifted so that it is centered at $s$. Finally, let $\mathbf{Y}_t = X_{U_t}$ and $\mathbf{Y}_t(s) = X_{U_t(s)}$. The dimension of $U_t$ varies for different $t$, depending on whether it is a corner or a middle pixel, but is always between $w(w-1)$ and $2w(w-1)$. So we will need a kernel $W^{(p)}$ for each $p$, $w(w-1) \leq p \leq 2w(w-1)$. Given the observed texture $\{X_t : t \in [1,T_1] \times [1,T_2]\}$, the algorithm to synthesize $X^*$ can be written as follows:

THE SPIRAL RESAMPLING ALGORITHM.

1. Select a random starting value for $\{X_t^* : \|t\|_\infty < w\}$, the central $(2w-1) \times (2w-1)$ square, uniformly from the observed field $X_t$.
2. Suppose $X_s^*$ have been generated for all $s \prec t$. To generate the next value $X_t^*$, let $N$ be a discrete random variable with probability mass function

$$P(N = s) = \frac{1}{Z} W_b^{(p)}(\mathbf{Y}_t^* - \mathbf{Y}_t(s)),$$

where $Z = \sum_s W_b^{(p)}(\mathbf{Y}_t^* - \mathbf{Y}_t(s))$ is a normalizing constant, $p = |\mathbf{Y}_t^*|$ is the size of the conditioning neighborhood for $X_t^*$, and $s$ ranges over all values in $[1,T_1] \times [1,T_2]$ such that $U_t(s) \in [1,T_1] \times [1,T_2]$.
3. Let $X_t^* = X_N$.

3.3. *Comparisons with the original algorithm of Efros and Leung and selection of tuning parameters.* In this section we investigate the effects of different orderings of the synthesized pixels, different resampling weights and tuning parameters. Only selected comparisons are shown here for the obvious reasons of space limitations; the conclusions drawn are based on a thorough simulation study comparing all variants on a larger number of images. For any particular comparison, all other parameters are held fixed at their optimal values.

Figure 5 shows the effects of changing the order, and also compares uniform versus kernel weights. The spiral order of the original Efros and Leung algorithm [Figure 5(b)] does appear to produce better results than the MMM version [Figure 5(e)], at least for the first texture (for the second mesh texture, all results are very similar). However, we claim that the difference is mainly due not to the spiral versus raster ordering, but to the fact that the spiral conditioning neighborhood contains twice as many close neighbors as the MMM "corner" neighborhood. To illustrate, we also generated textures in raster order [Figure 5(c), 5(d)] but conditioning on the full half-square above $t$, $U_t = \{u : \max(1, t_1 - w + 1) \leq u_1 \leq t_1, \max(1, t_2 - w + 1) \leq u_2 \leq t_2 + w - 1, u \neq t\}$, a version we will refer to as rectangular (as opposed to spiral and corner). This is not a MMM, and it generates results very similar to the original algorithm [Figure 5(b), 5(c)]. The remaining slight differences are probably due to the fact that spiral places the seed in the middle whereas the raster order starts from the corner, so error propagation is worse for a raster image of the same size. We also note that the speeded-up version of Efros and Leung [26] works in raster order with no problems, and so does the patch-based version in [8].

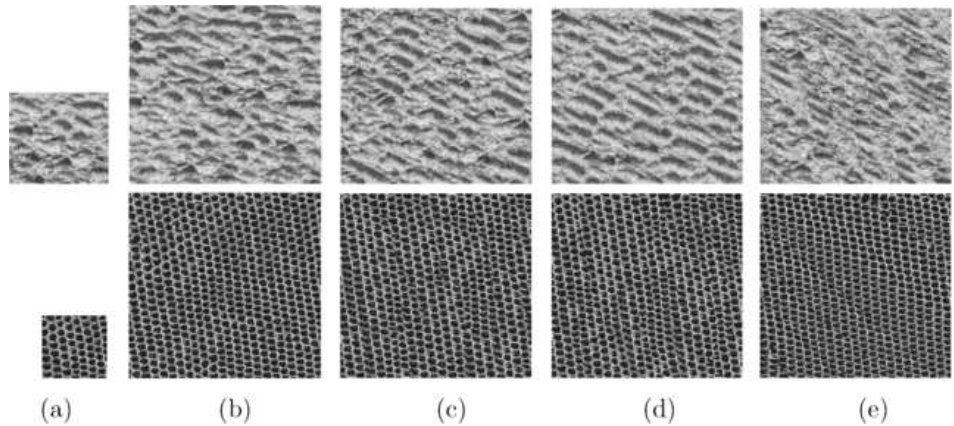

Fig. 5. *Comparing different orderings and uniform and kernel weights. (*a*) Original textures (sizes $151 \times 138$ and $81 \times 78$ pixels); (*b*) Efros and Leung result (spiral scheme with uniform weights); (*c*) rectangular scheme with uniform weights; (*d*) rectangular scheme with kernel weights; (*e*) MMM (corner scheme) with kernel weights. The window sizes are $w = 27$ and $w = 23$, respectively; the bandwidths are $\varepsilon = 0.1$ and $b = 0.01$.*

Using the uniform versus kernel weights [Figure 5(c), 5(d)] does not produce any detectable differences when $\varepsilon$ and $b$ are carefully chosen. The effect of varying $b$ is shown in Figure 6 and, predictably, increasing $b$ leads to the synthesized texture looking more "stochastic" and eventually becoming like white noise. Increasing $\varepsilon$ with uniform weights has the same effect. To make bandwidths choices more universal, all images are scaled to have grayscale values ranging from 0 to 1.

The effect of the window size has been shown in Figure 3, and remains the same for all versions. Other things being equal, larger window sizes tend to produce better results; however, they also make the computation costlier and reduce the effective sample size of the original image.

Finally, to keep things in perspective we note that all the different versions of the nonparametric resampling scheme [spiral with uniform weights vs. corner with kernel weights shown in Figure 7(b), 7(c)] are close to each other and quite good when compared to other texture synthesis methods, such as De Bonet [5] [Figure 7(d)] and Heeger and Bergen [10] [Figure 7(e)].

In the next section, we show that, subject to certain mixing and regularity conditions, both MMM and spiral schemes with kernel weights reproduce the joint distribution of a pixel in a patch consistently.

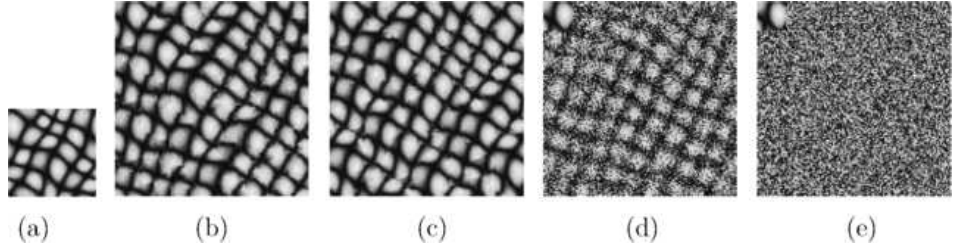

FIG. 6. *Kernel bandwidth effect. (a) The original texture (73 × 71 pixels); (b) $b = 0.007$ (smallest allowed by machine precision); (c) $b = 0.01$; (d) $b = 0.1$; (e) $b = 1$. All results generated with rectangular scheme, $w = 37$.*

**4. Consistency results.** We start by showing consistency of the Markov mesh model algorithm. For simplicity of notation, we ignore the side effects of truncation and show that the distribution of pixel value $X_t$ given its full $w$-by-$w$ neighborhood $\mathbf{Y}_t$ (for $t$ not in the first $w$ rows or columns) converges to the truth. We will then show that the same argument applies to truncated neighborhoods, and in fact to a neighborhood of an arbitrary shape, as long as the resampling scheme is matching it to neighborhoods of the same shape in the observed image. Finally, we generalize the consistency results to the original spiral ordering of the Efros and Leung algorithm.

4.1. *Assumptions.* Let us introduce the following notation: let $F_{\mathbf{Y}}(\mathbf{y}) = P(\mathbf{Y}_t \leq \mathbf{y})$ be the cumulative distribution function of $\mathbf{Y}_t$ and let $F_{X|\mathbf{Y}}(x|\mathbf{y}) = P(X_t \leq x|\mathbf{Y}_t = \mathbf{y})$ be the conditional distribution function of $X_t$ given $\mathbf{Y}_t$. Let $I = [w, T_1] \times [w, T_2]$ be the set of all pixels that admit a full conditional neighborhood. We assume that the size of the observed texture increases, that is, $T = \min(T_1, T_2) \to \infty$. We also make the following fairly technical regularity and mixing assumptions, which, however, are not unreasonable for real textures (see the discussion in Section 5). The assumption of com-pact support (A2) is automatically satisfied for images since the number of grayscale or color values used is finite.

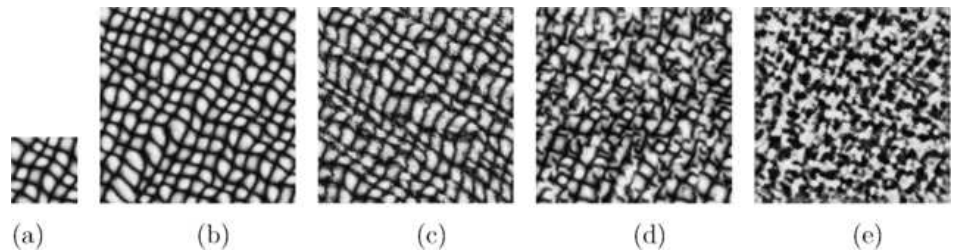

FIG. 7. *Nonparametric resampling compared to other methods. (a) The original texture (73 × 71 pixels); (b) Efros and Leung algorithm; (c) MMM with kernel weights ($w = 37$, $b = 0.01$); (d) De Bonet algorithm; (e) Heeger and Bergen algorithm.*

(A1) The random field $X_t$ is strictly stationary and mixing in the following sense: define mixing coefficients

$$\alpha_X(k,u,v) = \sup\{|P(AB) - P(A)P(B)| : A \in \sigma(X_E), B \in \sigma(X_F),$$
$$E, F \subset I, d(E,F) \geq k, |E| \leq u, |F| \leq v\},$$

where $d(E,F) = \inf\{\|x-y\|_\infty : x \in E, y \in F\}$ is the distance between index sets $E$ and $F$. The field $X_t$ is called $\alpha$-mixing if for all $u$ and $v$ $\alpha_X(k,u,v) \to 0$ as $k \to \infty$.

We make a more precise assumption about the rate at which the mixing coefficients go to 0: there exist $\varepsilon > 0$, $\tau > 2$ such that for all integers $u, v \geq 2, u+v \leq c$, where $c$ is the smallest even integer such that $c \geq \tau$,

$$\sum_{k=1}^{\infty} (k+1)^{d(c-u+1)-1} [\alpha_X(k,u,v)]^{\varepsilon/(c+\varepsilon)} < \infty.$$

Here $d$ is the dimension of the index set $I \subset \mathbb{Z}^d$, in our case $d = 2$.

(A2) $F_\mathbf{Y}$ and $F_{X|\mathbf{Y}}$ have bounded densities with respect to the Lebesgue measure, $f_\mathbf{Y}$ and $f_{X|\mathbf{Y}}$, respectively. Moreover, $X_t$ has compact support $S$, and $f_{X|\mathbf{Y}}(\cdot|\mathbf{y}) > 0$ for all $\mathbf{y} \in S^p$.

(A3) For any $\mathbf{y}_1, \mathbf{y}_2 \in \mathbb{R}^p$, any $x \in \mathbb{R} \cup \{\infty\}$,

$$\left| \int_{-\infty}^{x} f_{X,\mathbf{Y}}(z, \mathbf{y}_1) \, dz - \int_{-\infty}^{x} f_{X,\mathbf{Y}}(z, \mathbf{y}_2) \, dz \right| \leq L \|\mathbf{y}_1 - \mathbf{y}_2\|,$$

where $f_{X,\mathbf{Y}} = f_{X|\mathbf{Y}} f_\mathbf{Y}$.

(A4) The kernel $W$ on $\mathbb{R}^p$ is bounded, first-order Lipschitz continuous, symmetric, positive everywhere on $\mathbb{R}^p$, $\int uW(u)\,du = 0$, and $\int \|u\|W(u)\,du < \infty$. When $T \to \infty$, the kernel bandwidth $b = O(T^{-\delta})$, with $\delta > 0$ chosen so that $\delta < (\tau - 2)/2p(p+1+\tau)$.

4.2. *Consistency of the MMM algorithm.* Let $F^*_{X_t|\mathbf{Y}_t}(x|\mathbf{y}) = P(X^*_t \leq x | \mathbf{Y}^*_t = \mathbf{y})$ be the conditional distribution function of the synthesized $X^*_t$ given its neighborhood $\mathbf{Y}^*_t = \mathbf{y}$ and let $F^*_{X_t, \mathbf{Y}_t}(x, \mathbf{y}) = P(X^*_t \leq x, \mathbf{Y}^*_t \leq \mathbf{y})$ be the joint distribution function of $X^*_t$ and $\mathbf{Y}^*_t$, that is, the joint distribution of pixels in a $w \times w$ window. This is the distribution of interest because, at least for some suitably chosen $w$, it determines the human perception of texture, as discussed in the Introduction. Therefore one may argue that if this joint distribution is estimated correctly, then the synthesized texture will appear similar to the original. Our main result is the following theorem.

THEOREM 1. *Under assumptions (A1)–(A4), the joint distribution of $X_t^*$ and $\mathbf{Y}_t^*$ (the joint distribution of pixels in a $w \times w$ window) is estimated consistently for all $t \in [w, \infty)^2$:*

$$\sup_{x \in \mathbb{R}} \sup_{\mathbf{y} \in S^p} |F^*_{X_t, \mathbf{Y}_t}(x, \mathbf{y}) - F_{X, \mathbf{Y}}(x, \mathbf{y})| \to 0 \qquad \text{a.s. as } T \to \infty. \tag{1}$$

We also prove that the resampling scheme correctly approximates the conditional distribution of a pixel given its neighborhood.

THEOREM 2. *Under assumptions (A1)–(A4), the conditional distribution of $X_t^*$ given $\mathbf{Y}_t^*$ (the distribution of the right bottom corner pixel in a $w \times w$ window given the other $p = w^2 - 1$ pixels in the window) is estimated consistently for all $t \in [w, \infty)^2$:*

$$\sup_{x \in \mathbb{R}} \sup_{\mathbf{y} \in S^p} |F^*_{X_t | \mathbf{Y}_t}(x | \mathbf{y}) - F_{X | \mathbf{Y}}(x | \mathbf{y})| \to 0 \qquad \text{a.s. as } T \to \infty. \tag{2}$$

These theorems establish the consistency of the joint distribution of pixels in a $w \times w$ window. Inspection of the proof shows that the argument does not depend on the shape of $\mathbf{Y}_t$. All it requires is that the number of observed $\mathbf{Y}_t$ goes to infinity, so that there are many matches to sample from. It also does not depend on the particular order in which the pixels are synthesized, because the argument is for a single given pixel in the synthesized texture as the size of the *observed* texture grows. If in the beginning the seed is chosen uniformly from the original, we start from a set of pixels whose joint distribution is consistent, and add pixels one by one in such a way that the joint distribution of the $w \times w$ window with that pixel in the bottom right corner remains consistent. Thus the joint distribution of every $w \times w$ window throughout the synthesized texture is estimated consistently. This is the main result we were interested in, since it suggests that the synthesized texture will appear similar to the original.

4.3. *Consistency of the spiral resampling algorithm.* It is clear from the proofs of Theorems 1 and 2 that as long as we assume the mixing assumption (A1) holds, all kernels $W^{(p)}$ satisfy assumption (A4), and for all shapes of $\mathbf{Y}$ the distributions of $X$ and $\mathbf{Y}$ satisfy assumptions (A2) and (A3), we will obtain the same consistency result: the conditional distribution of $X_t^*$ given whatever pixels $\mathbf{Y}_t^*$ we conditioned on to fill it in will converge to the true distribution of $X$ given $\mathbf{Y}$ as the size of the observed texture goes to infinity. Similarly, the joint distribution of $X_t^*$ and $\mathbf{Y}_t^*$ will converge to the truth, at every pixel location in the synthesized texture. However, in this case the shape of the neighborhood depends on location and is constantly changing according to the spiral ordering. Therefore it is not clear whether we can obtain a consistent estimate of the joint distribution in a square window over the whole synthesized texture.

It would be natural to expect that, since we always condition on a $(2w-1) \times (2w-1)$ window, or at least on what we can see from it, we will in the end get the joint distribution estimate in that window consistently. This was the motivation for the Efros and Leung algorithm, but it is in fact not true.

Consider a simple counterexample for $w = 2$ shown in Figure 8. The first pixel to be filled in after the seed is $X_1^*$, and the part of the $3 \times 3$ window around it that we have so far is $S_1^*$. We assume that the seed was sampled uniformly from the original image, so the sampling distribution $P(S_1^*) \to P(S_1)$. Since the conditional distribution of $X_1^*$ given $S_1^*$ is estimated consistently, we have

$$P(X_1^*, S_1^*) = P(X_1^*|S_1^*)P(S_1^*) \to P(X_1|S_1)P(S_1) = P(X_1, S_1),$$

so the joint distribution of $X_1^*$ and $S_1^*$ is estimated consistently. However, this tells us nothing about the joint distribution of $X_1^*$ and $S_2^*$. In fact, by construction of the synthesis algorithm we have

$$\begin{aligned} P(X_1^*, S_1^*, S_2^*) &= P(X_1^*, S_2^*|S_1^*)P(S_1^*) = P(X_1^*|S_1^*)P(S_2^*|S_1^*)P(S_1^*) \\ &\to P(X_1|S_1)P(S_2|S_1)P(S_1) = P(X_1|S_1)P(S_1, S_2) \\ &\neq P(X_1|S_1, S_2)P(S_1, S_2) = P(X_1, S_1, S_2). \end{aligned}$$

In other words, the synthesis algorithm makes $X_1^*$ and $S_2^*$ independent given $S_1^*$, a property that the true distribution does not have in general. Therefore the estimate of the joint distribution in a $(2w-1) \times (2w-1)$ window cannot be consistent. However, we may still get consistency in a smaller window, and in fact this example suggests that in order to get consistency for a window of size $w \times w$, one must condition on a bigger window that contains all $w \times w$ windows or their parts which cover the pixel being synthesized. The size of this bigger window must be exactly $(2w-1) \times (2w-1)$ in order to cover all $w \times w$ windows containing the pixel at its center. Then each added pixel will fit in correctly with all $w \times w$ windows that contain it, and by induction the joint distribution in all $w \times w$ windows throughout the synthesized texture will be estimated consistently.

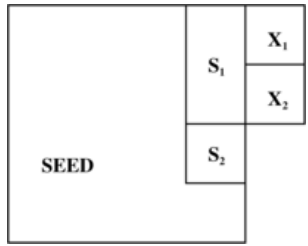

FIG. 8. *Counterexample to consistency of the original algorithm: in the bootstrapped version, $X_1$ and $S_2$ are independent given $S_1$.*

To formalize this claim, let

$$V_t = \{s : \|t - s\|_\infty < w, s \neq t\}$$

be the $(2w-1) \times (2w-1)$ window centered at $t$, and let $\mathbf{V}_t = X_{V_t}$ be the pixel intensities in that window. Let $Q_t$ be the pixels in the $w \times w$ window located at $[t - w + 1, t]^2$, and let $\mathbf{Q}_t = X_{Q_t}$ be the corresponding vector of pixel intensities. Let

$$F^*_{\mathbf{Q}_t}(\mathbf{q}) = P(\mathbf{Q}^*_t \leq \mathbf{q})$$

be the cumulative distribution function of $\mathbf{Q}_t$ in the texture synthesized by the spiral-order algorithm, and let $F_\mathbf{Q}$ be the true cumulative distribution function of the $w \times w$ window in the original texture.

THEOREM 3. *Suppose we observe $\{X_t : t \in [1, T_1] \times [1, T_2]\}$. Let $T = \min(T_1, T_2)$. If the field $X_t$ satisfies assumption (A1), the distributions of $X$, $\mathbf{V}$, $X|\mathbf{V}$ and $(X, \mathbf{V})$ satisfy assumptions (A2) and (A3), and the kernels $W^{(p)}$ for all $p = w(w-1), \ldots, 2w(w-1)$ satisfy assumption (A4), then the distribution of $\mathbf{Q}_t$ is estimated consistently for all $t \in \mathbb{Z}^2$:*

(3) $$\sup_{q \in S^{w^2}} |F^*_{\mathbf{Q}_t}(\mathbf{q}) - F_\mathbf{Q}(\mathbf{q})| \to 0 \qquad a.s. \text{ as } T \to \infty.$$

This shows that the original algorithm of Efros and Leung also provides consistent estimates of the joint distribution in a $w \times w$ window, which may be an explanation for its perceptually good results, although the window is smaller than what the authors intended. One can similarly show that conditioning on a rectangular upper half-window when synthesizing in raster order (what we called the "rectangular scheme" in experiments) produces a consistent estimate of the distribution in a $w \times w$ window, but not in the full $(2w-1) \times (2w-1)$ window.

**5. Discussion.** The main contributions of this paper are the formal probabilistic framework for the nonparametric sampling algorithm of Efros and Leung [9] and the proof of its consistency. In particular, the fact that the joint distribution of pixels in a window of specified size is estimated consistently may explain the perceptually good results of the algorithm. This joint distribution is important for texture perception both from the Julesz school's point of view ($k$th-order statistics) and from the multichannel frequency analysis perspective, since the joint distribution of pixels in the filters' support window determines the joint distribution of filter responses.

The proof of consistency requires a number of conditions which may look complicated, but are in fact perfectly plausible for most real textures. The mixing assumption will not hold for purely periodic patterns, but as long as there is some stochastic variation, it becomes a natural description of texture. As for the density assumptions, the grayscale pixel values are discrete, but there are sufficiently many of them to make the smoothness conditions plausible. The assumption of compact support is automatically satisfied for images since the number of grayscale or color values used is always finite, but may be a more substantial limitation for applications to other random fields.

Our goal here was to provide asymptotic justification of the basic nonparametric resampling principle rather than of any particular implementation of it. For instance, the resampling weights we use are slightly different from the ones originally used by Efros and Leung [9]. We feel that modifying the proof to accommodate details of a particular implementation, either the original or one of the several follow-up versions, is possible but unnecessary, since this is intended as a justification of their common underlying principle.

One issue that we did not address is determining the correct window size to use in the resampling algorithm in order to obtain perceptually good results. The asymptotics guarantee only that the distribution in a certain size window is estimated consistently. They say nothing about whether the distribution of the whole synthesized texture is consistent, unless we are willing to assume that the true distribution is a Markov random field with the neighborhood of the same or smaller size as our chosen window. Experimentally it appears that the window big enough to contain the largest texture "element" (determined by the user) works for the resampling algorithm. Automatically determining the correct window size for nonparametric resampling algorithms, and the scale of a given texture in general, is an open problem. In classic MBB, the size of the block can be chosen to optimize the bias and variance of the estimator of the mean; this cannot be applied here as the goal is to reproduce the joint distribution of pixels rather than to estimate the mean. One could use cross-validation, that is, compare the synthesized texture to the original for several window sizes using a texture similarity measure (see, e.g., [16]), and pick the window size that maximizes this similarity. This approach is somewhat computationally expensive, and there is no guarantee that the similarity measures used for classification and segmentation will be adequate for human perception.

Another tuning parameter set by the user is the kernel bandwidth $b$ or $\varepsilon$. In the algorithm implementations they were determined empirically and held constant for all the textures, so it only needed to be done once. Methods for bandwidth selection used in density estimation could be applied here, although one does not expect drastic practical improvements. Another possibility is to select both the window size and the kernel bandwidth by cross-validation, which may yet become preferred over the standard computer vision practice of user-selected parameters as computing becomes faster and similarity measures get fine-tuned to mimic human perception more closely.

A natural question to ask is what else this algorithm could be useful for beyond texture synthesis. We believe it will do well for the usual bootstrap task of estimating the mean and variance of a random field; establishing its rates of convergence and comparing them to, for example, MBB could be a direction for future work. A drawback from the computer vision point of view is that this type of bootstrap cannot be used to perform texture classification or recognition, since no generative model is fit to the data. It can, however, be used for estimating various texture parameters, such as the texture scale, via cross-validation as described above; these parameters can in turn be useful for classification and other higher-level tasks.

## APPENDIX: PROOFS

Everywhere in the proofs we suppress the dependence on $t$ in $F^*_{X_t,\mathbf{Y}_t}$ and $F^*_{X_t|\mathbf{Y}_t}$ to avoid clutter, and remind the reader all theorems hold for all appropriate $t \in \mathbb{Z}^2$. Before proceeding to the proofs of our results, we state a moment inequality for mixing random fields which we will need below. The proof of this inequality and many other useful ones can be found in [7].

LEMMA A.1 (Moment inequality). *Let $F_t$ be a real-valued random field indexed by $I \subset \mathbb{Z}^d$ satisfying conditions (A1). If $EF_t = 0$, $F_t \in L^{\tau+\varepsilon}$ and $\tau \geq 2$, then there is a constant $C$ depending only on $\tau$ and mixing coefficients of $F_t$ such that*

$$E\left|\sum_{t \in I} F_t\right|^\tau \leq C \max(L(\tau,\varepsilon), L(2,\varepsilon)^{\tau/2}),$$

*where*

$$L(\mu,\varepsilon) = \sum_{t \in I}(E|F_t|^{\mu+\varepsilon})^{\mu/(\mu+\varepsilon)} = \sum_{t \in I} \|F_t\|^\mu_{\mu+\varepsilon}.$$

We will start from the proof of Theorem 2 (consistency of conditional distributions for the MMM algorithm). Note that, for any $x \in \mathbb{R}$, $\mathbf{y} \in S^p$,

$$F^*_{X|\mathbf{Y}}(x|\mathbf{y}) = \sum_{t \in I} \mathbf{1}_{(-\infty,x]}(X_t)W_b(\mathbf{y}-\mathbf{Y}_t) \Big/ \sum_{s \in I} W_b(\mathbf{y}-\mathbf{Y}_s),$$

$$F_{X|\mathbf{Y}}(x|\mathbf{y}) = \int \mathbf{1}_{(-\infty,x]}(z) f_{X|\mathbf{Y}}(z|\mathbf{y})\, dz.$$

We first prove the following unconditional result.

LEMMA A.2. *Under assumptions (A1)–(A4), for any $x \in \mathbb{R}$*

$$\sup_{\mathbf{y} \in S^p} \left| \frac{1}{T} \sum_{t \in I} \mathbf{1}_{(-\infty, x]}(X_t) W_b(\mathbf{y} - \mathbf{Y}_t) - \int \mathbf{1}_{(-\infty, x]}(z) f_{X, \mathbf{Y}}(z, \mathbf{y}) \, dz \right| \to 0 \quad \text{a.s.}$$

From this lemma, we can immediately get a useful corollary. Let $f^*_\mathbf{Y}(\mathbf{y}) = \frac{1}{T} \sum_{t \in I} W_b(\mathbf{y} - \mathbf{Y}_t)$. Then setting $x = \infty$ in Lemma A.2 we get:

COROLLARY A.1. *Under the conditions of Theorem 2,*

$$\sup_{\mathbf{y} \in S^p} |f^*_\mathbf{Y}(\mathbf{y}) - f_\mathbf{Y}(\mathbf{y})| \to 0 \quad \text{a.s.}$$

PROOF OF LEMMA A.2. Let us introduce the notation

$$r^*_T(\mathbf{y}) = \frac{1}{T} \sum_{t \in I} \mathbf{1}_{(-\infty, x]}(X_t) W_b(\mathbf{y} - \mathbf{Y}_t),$$

$$r(\mathbf{y}) = \int \mathbf{1}_{(-\infty, x]}(z) f_{X, \mathbf{Y}}(z, \mathbf{y}) \, dz.$$

The lemma will be proved by showing

(4) $$\sup_{\mathbf{y} \in S^p} |Er^*_T(\mathbf{y}) - r(\mathbf{y})| \to 0 \quad \text{a.s.}$$

and

(5) $$\sup_{\mathbf{y} \in S^p} |r^*_T(\mathbf{y}) - Er^*_T(\mathbf{y})| \to 0 \quad \text{a.s.}$$

First let us compute

$$Er^*_T(\mathbf{y}) = \int_{-\infty}^{x} \int_{\mathbf{u} \in R^p} f_{X, \mathbf{Y}}(z, \mathbf{u}) W_b(\mathbf{y} - \mathbf{u}) \, d\mathbf{u} \, dz$$

$$= b^{-p} \int_{-\infty}^{x} \int_{\mathbf{u} \in R^p} f_{X, \mathbf{Y}}(z, \mathbf{u}) W\left(\frac{\mathbf{u} - \mathbf{y}}{b}\right) d\mathbf{u} \, dz$$

$$= \int_{-\infty}^{x} \int_{\mathbf{v} \in R^p} f_{X, \mathbf{Y}}(z, b\mathbf{v} + \mathbf{y}) W(\mathbf{v}) \, d\mathbf{v} \, dz,$$

where $\mathbf{v} = (\mathbf{u} - \mathbf{y})/b$. Also note that, since $W$ is a kernel,

$$r(\mathbf{y}) = \int_{-\infty}^{x} f_{X, \mathbf{Y}}(z, \mathbf{y}) \, dz = \int_{-\infty}^{x} \int_{\mathbf{v} \in R^p} f_{X, \mathbf{Y}}(z, \mathbf{y}) W(\mathbf{v}) \, d\mathbf{v} \, dz.$$

Now (4) becomes

$$\sup_{\mathbf{y} \in S^p} |Er^*_T(\mathbf{y}) - r(\mathbf{y})|$$

$$= \sup_{\mathbf{y} \in S^p} \left| \int_{-\infty}^{x} \int_{\mathbf{v} \in R^p} W(\mathbf{v})(f_{X, \mathbf{Y}}(z, b\mathbf{v} + \mathbf{y}) - f_{X, \mathbf{Y}}(z, \mathbf{y})) \, d\mathbf{v} \, dz \right|$$

$$\leq \sup_{\mathbf{y}\in S^p} \int_{\mathbf{v}\in R^p} W(\mathbf{v}) \left| \int_{-\infty}^{x} f_{X,\mathbf{Y}}(z, b\mathbf{v}+\mathbf{y})\,dz - \int_{-\infty}^{x} f_{X,\mathbf{Y}}(z, \mathbf{y})\,dz \right| d\mathbf{v}$$

$$\leq bL \int \|\mathbf{v}\| W(\mathbf{v})\,d\mathbf{v} = O(b) = O(T^{-\delta}) \to 0 \quad \text{a.s.},$$

where the last inequality follows from assumption (A3).

Equation (5) is the main part of the proof. Define

$$Z_{t,T}(\mathbf{y}) = \mathbf{1}_{(-\infty,x]}(X_t)W_b(\mathbf{y}-\mathbf{Y}_t) - E\mathbf{1}_{(-\infty,x]}(X_t)W_b(\mathbf{y}-\mathbf{Y}_t).$$

Note that $EZ_{t,T} = 0$, and claim (5) is that

$$\sup_{\mathbf{y}\in S^p} \left| \frac{1}{T} \sum_{t\in I} Z_{t,T}(\mathbf{y}) \right| \to 0 \quad \text{a.s.}$$

Recall that $X_t$ has compact support $S$. Therefore we can cover $S^p$ with $N_T$ cubes $I_{i,T}$ with centers $y_i$ and sides $L_T$. Then

$$\sup_{\mathbf{y}\in S^p} \left| \frac{1}{T} \sum_{t\in I} Z_{t,T}(\mathbf{y}) \right|$$

$$= \max_{1\leq i \leq N_T} \sup_{\mathbf{y}\in S^p \cap I_{i,T}} \left| \frac{1}{T} \sum_{t\in I} Z_{t,T}(\mathbf{y}) \right|$$

$$\leq \max_{1\leq i \leq N_T} \left| \frac{1}{T} \sum_{t\in I} Z_{t,T}(\mathbf{y}_i) \right|$$

$$+ \max_{1\leq i \leq N_T} \sup_{\mathbf{y}\in S^p \cap I_{i,T}} \left| \frac{1}{T} \sum_{t\in I} (Z_{t,T}(\mathbf{y}) - Z_{t,T}(\mathbf{y}_i)) \right|$$

$$= \text{I} + \text{II}.$$

First let us deal with term II:

$$\text{II} \leq \max_{1\leq i \leq N_T} \sup_{\mathbf{y}\in S^p \cap I_{i,T}} \frac{1}{T} \left[ \sum_{t\in I} |W_b(\mathbf{y}-\mathbf{Y}_t) - W_b(\mathbf{y}_i-\mathbf{Y}_t)| \right.$$

$$\left. + \sum_{t\in I} E|W_b(\mathbf{y}-\mathbf{Y}_t) - W_b(\mathbf{y}_i-\mathbf{Y}_t)| \right]$$

$$\leq C_1 \max_{1\leq i \leq N_T} \sup_{\mathbf{y}\in S^p \cap I_{i,T}} b^{-p} \left\| \frac{\mathbf{y}-\mathbf{y}_i}{b} \right\|$$

$$\leq C_1 b^{-p-1} L_T.$$

The last line follows from the Lipschitz assumption on the kernel (A4).

If we let the side of the cubes $L_T = \varepsilon b^{p+1} = O(T^{-\delta(p+1)})$, then term II is bounded above by $\varepsilon$. Note that the number of cubes $N_T = O(1/L_T^p) = O(T^{\delta p(p+1)})$.

We will use the Borel–Cantelli lemma to show that term I goes to 0. By elementary inequalities,

$$P_T(\mathrm{I} > \varepsilon) \leq \sum_{i=1}^{N_T} P\left(\left|\frac{1}{T}\sum_{t \in I} Z_{t,T}(\mathbf{y}_i)\right| > \varepsilon\right)$$

(6)
$$\leq N_T \max_{1 \leq i \leq N_T} P\left(\left|\frac{1}{T}\sum_{t \in I} Z_{t,T}(\mathbf{y}_i)\right| > \varepsilon\right)$$

$$\leq N_T \max_{1 \leq i \leq N_T} \frac{E|\sum_t Z_{t,T}(\mathbf{y}_i)|^\tau}{\varepsilon^\tau T^\tau}.$$

To bound the last term, we apply the moment inequality of Lemma A.1 to random variables,

$$F_t = Z_{t,T} = \mathbf{1}_{(-\infty,x]}(X_t) W_b(\mathbf{y} - bY_t) - E\mathbf{1}_{(-\infty,x]}(X_t) W_b(\mathbf{y} - bY_t),$$

so we need to check that they satisfy assumption (A1). First note that

$$F_t = f(X_t, \mathbf{Y}_t) = f(X_{[t-w+1,t]^2}).$$

Since the definition of mixing coefficients only depends on $\sigma$-algebras generated by $F_t$, we may instead consider larger $\sigma$-algebras generated by $(X_t, \mathbf{Y}_t)$—the mixing coefficients of $F_t$ can only be smaller than those of $(X_t, \mathbf{Y}_t)$. Therefore,

$$\alpha_F(k, u, v) \leq \sup\{|P(AB) - P(A)P(B)| : A \in \sigma((X, \mathbf{Y})_E), B \in \sigma((X, \mathbf{Y})_F),$$
$$E, F \subset I, d(E, F) \geq k, |E| \leq u, |F| \leq v\}.$$

Now notice that $(X, \mathbf{Y})_E = X_{E'}$ where $E' = \{t + a \in I : t \in E, a \in [-w, 0] \times [-w, 0]\}$, that is, all points in $E$ and everything in a $w$-by-$w$ square to the left and above them. Similarly, let $F' = \{t + a \in I : t \in F, a \in [-w, 0] \times [-w, 0]\}$. If $|E| \leq u$ and $|F| \leq v$, then $|E'| \leq w^2 u$ and $|F'| \leq w^2 v$. Also, if the distance between $E$ and $F$, $d(E, F) \geq k$, then $d(E', F') \geq k - w$. Therefore,

$$\alpha_F(k, u, v) \leq \alpha_X(k - w, w^2 u, w^2 v),$$

and assumptions (A1) are clearly satisfied for the field $F_t$, possibly with different constants.

Now we can apply Lemma A.1 to identically distributed mixing variables $Z_{t,T}$. Note that since the kernel $W$ is bounded, $|Z_{t,T}(\mathbf{y}_i)| \leq \tilde{M} b^{-p} \leq M T^{\delta p}$. Also note

$$L(\tau, \varepsilon) = T(E|Z_{t,T}|^{\tau+\varepsilon})^{\tau/(\tau+\varepsilon)} \leq T(MT^{\delta p})^\tau = M^\tau T^{\delta p \tau + 1}$$

and
$$L(2,\varepsilon)^{\tau/2} \leq (M^2 T^{2\delta p+1})^{\tau/2} = M^\tau T^{\delta p\tau + \tau/2},$$
and since we take $\tau \geq 2$
$$E\left|\sum_{t\in I} Z_{t,T}\right| \leq C\max(M^\tau T^{\delta p\tau+1}, M^\tau T^{\delta p\tau+\tau/2}) \leq \tilde{C}T^{\delta p\tau+\tau/2}.$$

Now we can go back to (6) and plug in the moment bound. We get
$$P_T(\mathrm{I} > \varepsilon) \leq N_T \max_{1\leq i \leq N_T} \frac{C(\tau)T^{\delta p\tau + \tau/2}}{\varepsilon^\tau T^\tau} = O(T^{\delta p(p+1) + \delta p\tau - \tau/2}).$$

Since we assumed $\tau > 2$ and $\delta < (\tau - 2)/2p(p+1+\tau)$, we can always choose $\delta$ such that $\delta p(p+1) + \tau(\delta p - 1/2) < -1$. Therefore $\sum_T P_T(I > \varepsilon) < \infty$ and by the Borel–Cantelli lemma $I \to 0$ a.s. This concludes the proof of Lemma A.2. □

PROOF OF THEOREM 2. Note that
$$F^*_{X|\mathbf{Y}}(x|\mathbf{y}) = \frac{\sum_{t\in I}\mathbf{1}_{(-\infty,x]}(X_t)W_b(\mathbf{y}-\mathbf{Y}_t)}{\sum_{t\in I} W_b(\mathbf{y}-\mathbf{Y}_t)} = \frac{r^*_T(\mathbf{y})}{f^*_\mathbf{Y}(\mathbf{y})}$$
and that $r(\mathbf{y}) = F_{X|\mathbf{Y}}(x|\mathbf{y})f_\mathbf{Y}(\mathbf{y})$. Then the expression under sup in (2) becomes
$$|F^*_{X|\mathbf{Y}}(x|\mathbf{y}) - F_{X|\mathbf{Y}}(x|\mathbf{y})|$$
$$= \frac{1}{f^*_\mathbf{Y}(\mathbf{y})}|r^*_T(\mathbf{y}) - r(\mathbf{y}) + r(\mathbf{y}) - F_{X|\mathbf{Y}}(x|\mathbf{y})f^*_\mathbf{Y}(\mathbf{y})|$$
$$\leq \frac{1}{f^*_\mathbf{Y}(\mathbf{y})}(|r^*_T(\mathbf{y}) - r(\mathbf{y})| + F_{X|\mathbf{Y}}(x|\mathbf{y})|f^*_\mathbf{Y}(\mathbf{y}) - f_\mathbf{Y}(\mathbf{y})|).$$

From Lemma A.2, Corollary A.1 and assumption (A2) it follows that
$$\sup_{\mathbf{y}\in S^p}|F^*_{X|\mathbf{Y}}(x|\mathbf{y}) - F_{X|\mathbf{Y}}(x|\mathbf{y})| \to 0 \quad \text{a.s.}$$

To establish the uniform convergence over all $x \in \mathbb{R}$, we use the argument of the Glivenko–Cantelli theorem: for each $x$, $F^*_{X|\mathbf{Y}}(x|\mathbf{y}) \to F_{X|\mathbf{Y}}(x|\mathbf{y})$ a.s. by the ergodic theorem. Since $F^*_{X|\mathbf{Y}}$ are nondecreasing, and $F_{X|\mathbf{Y}}$ is bounded and continuous, it follows that convergence is uniform over all $x \in \mathbb{R}$. □

PROOF OF THEOREM 1. First note that
$$F^*_{X,\mathbf{Y}}(x,\mathbf{y}) = P(X^*_t \leq x, \mathbf{Y}^*_t \leq \mathbf{y}) = \int_{\mathbf{u}\leq \mathbf{y}}\int_{z\leq x} f^*_{X,\mathbf{Y}}(z,\mathbf{u})\,dz\,d\mathbf{u}$$
$$= \int_{\mathbf{u}\leq \mathbf{y}}\int_{z\leq x} f^*_{X|\mathbf{Y}}(z|\mathbf{u})f^*_\mathbf{Y}(\mathbf{u})\,dz\,d\mathbf{u} = \int_{\mathbf{u}\leq \mathbf{y}} F^*_{X|\mathbf{Y}}(x|\mathbf{u})f^*_\mathbf{Y}(\mathbf{u})\,d\mathbf{u}.$$

Similarly,
$$F_{X,\mathbf{Y}}(x,\mathbf{y}) = \int_{\mathbf{u} \leq \mathbf{y}} F_{X|\mathbf{Y}}(x|\mathbf{u}) f_{\mathbf{Y}}(\mathbf{u}) \, d\mathbf{u}.$$

Recall that we assumed $f_\mathbf{Y}$ is bounded and has compact support $S^p$. Then we have the bound (indexes are omitted for clarity)

$$|F^*(x,\mathbf{y}) - F(x,\mathbf{y})|$$
$$\leq \int_{\mathbf{u} \leq \mathbf{y}} F^*(x|\mathbf{u})|f^*(\mathbf{u}) - f(\mathbf{u})| \, d\mathbf{u} + \int_{\mathbf{u} \leq \mathbf{y}} |F^*(x|\mathbf{u}) - F(x|\mathbf{u})| f(\mathbf{u}) \, d\mathbf{u}$$
$$\leq 1 \cdot \int |f^*(\mathbf{u}) - f(\mathbf{u})| \, d\mathbf{u} + M \int |F^*(x|\mathbf{u}) - F(x|\mathbf{u})| \, d\mathbf{u}$$
$$\leq |S^p| \sup_{\mathbf{u} \in S^p} |f^*(\mathbf{u}) - f(\mathbf{u})| + M|S^p| \sup_{\mathbf{u}} |F^*(x|\mathbf{u}) - F(x|\mathbf{u})|.$$

Taking the supremum over $x$ and $\mathbf{y}$, we get

$$\sup_{x \in \mathbb{R}} \sup_{\mathbf{y} \in S^p} |F^*(x,\mathbf{y}) - F(x,\mathbf{y})|$$
$$\leq C_1 \sup_{\mathbf{y} \in S^p} |f^*(\mathbf{y}) - f(\mathbf{y})| + C_2 \sup_{x \in \mathbb{R}} \sup_{\mathbf{y} \in S^p} |F^*(x|\mathbf{y}) - F(x|\mathbf{y})| \to 0 \quad \text{a.s.}$$

by Theorem 2 and Corollary A.1. $\square$

PROOF OF THEOREM 3. Let $\mathbf{Y}_t = \{X_s : s \in V_t, s \prec t\}$ be the vector of pixel values in $V_t$ that come before $X_t$ in the spiral ordering. The first thing we need to verify is that if $X$ and $\mathbf{V}$ [the full $(2w-1) \times (2w-1)$ window around $X$] satisfy assumptions (A2) and (A3), so will $X$ and $\mathbf{Y}$.

Assumption (A2) says that $\mathbf{V}$ and $X|\mathbf{V}$ have bounded, positive everywhere densities with compact support. This will obviously hold for $\mathbf{Y}$ and $X|\mathbf{Y}$ as well.

Now let us verify that assumption (A3) holds for $X$ and $\mathbf{Y}$. Write $\mathbf{V} = (\mathbf{Y}, \mathbf{U})$, where $\mathbf{U}$ is the part of $\mathbf{V}$ that comes after $X$. Then we can write

$$f_{X,\mathbf{Y}}(x,\mathbf{y}) = \int f_{X,\mathbf{Y},\mathbf{U}}(x,\mathbf{y},\mathbf{u}) \, d\mathbf{u}.$$

Let $\mathbf{v}_1 = (\mathbf{y}_1, \mathbf{u})$ and $\mathbf{v}_2 = (\mathbf{y}_2, \mathbf{u})$. Then

$$\left| \int_{-\infty}^{x} f_{X,\mathbf{Y}}(r,\mathbf{y}_1) \, dr - \int_{-\infty}^{x} f_{X,\mathbf{Y}}(r,\mathbf{y}_2) \, dr \right|$$
$$\leq \int_{\mathbf{u}} \left| \int_{-\infty}^{x} f_{X,\mathbf{Y},\mathbf{U}}(r,\mathbf{y}_1,\mathbf{u}) \, dr - \int_{-\infty}^{x} f_{X,\mathbf{Y},\mathbf{U}}(r,\mathbf{y}_2,\mathbf{u}) \, dr \right| d\mathbf{u}$$
$$\stackrel{(A3)}{\leq} |S|^{\dim(\mathbf{u})} L \|\mathbf{v}_1 - \mathbf{v}_2\| \leq \tilde{L} \|\mathbf{y}_1 - \mathbf{y}_2\|.$$

Now we can use Theorem 1 to conclude that

$$(7) \qquad \sup_x \sup_{\mathbf{y}} |F^*_{X,\mathbf{Y}}(x,\mathbf{y}) - F_{X,\mathbf{Y}}(x,\mathbf{y})| \to 0 \qquad \text{a.s. as } T \to \infty.$$

The consistency of the distribution estimate in a $w \times w$ window can now be proved by induction. Let $\mathcal{Q}_t$ be the set of all $w \times w$ windows and parts of such windows that are filled in by the time the pixel at location $t$ is synthesized. The induction hypothesis is that, for all $t \in \mathbb{Z}^2$, if $Q \in \mathcal{Q}_t$ and $\mathbf{Q} = X_Q$, then

$$(8) \qquad \sup_{\mathbf{q}} |F^*_{\mathbf{Q}}(\mathbf{q}) - F_{\mathbf{Q}}(\mathbf{q})| \to 0 \qquad \text{a.s. as } T \to \infty.$$

(1) For the first $(2w-1)^2$ locations $t_i$, $i = 0, \ldots, 4w(w-1)$, the windows in $\mathcal{Q}_t$ are sampled uniformly from the observed texture, and therefore $F^*_{\mathbf{Q}}$ is the empirical distribution function of the corresponding window in the texture sample. The claim (8) is true by the Glivenko–Cantelli theorem.

(2) Suppose (8) holds for all $s \prec t$. For all $Q \in \mathcal{Q}_t$ that do not include $t$, the claim holds since these sets also belong to $\mathcal{Q}_s$ for some $s \prec t$. For sets $Q$ that include $t$, we can write $\mathbf{Q}^* = (X^*_t, \mathbf{S}^*)$, where $\mathbf{S}^*$ is the vector of all pixel intensities in $Q$ other than $X^*_t$. Since the size of $Q$ is at most $w \times w$, by construction of our conditioning neighborhood all the pixels of $\mathbf{S}^*$ are included in $\mathbf{Y}^*_t$. By (7), the joint distribution of $X^*_t$ and $\mathbf{Y}^*_t$ converges to the truth, and therefore so does the joint distribution of $X^*_t$ and $S^*_t$. This establishes the induction hypothesis for all $Q \in \mathcal{Q}_t$. $\square$

**Acknowledgments.** Most of this work was done while the first author was a graduate student at University of California, Berkeley. The authors thank Professor Jitendra Malik and his former students Alyosha Efros and Jianbo Shi for helpful discussions, and Ed Ionides for a key suggestion in one of the proofs. The results of the Efros and Leung algorithm are reproduced with permission from the website of Alyosha Efros.

DEPARTMENT OF STATISTICS
UNIVERSITY OF MICHIGAN
ANN ARBOR, MICHIGAN 48109-1107
USA
E-MAIL: elevina@umich.edu

DEPARTMENT OF STATISTICS
UNIVERSITY OF CALIFORNIA, BERKELEY
BERKELEY, CALIFORNIA 94720-3860
USA
E-MAIL: bickel@stat.berkeley.edu